\documentclass[oneside, reqno] {amsart}

\usepackage{xypic}
\input xy
\xyoption{all}
\usepackage{epsfig}
\usepackage{amsthm}
\usepackage{amssymb}
\usepackage{amsmath}
\usepackage{amscd}
\usepackage{color}
\usepackage[T1]{fontenc}
\usepackage{upgreek}
\usepackage[font=scriptsize]{caption}
\usepackage{wrapfig}
\usepackage{graphicx}
\usepackage{subfigure}

\newcommand{\bg}{\begin{equation}}
\newcommand{\ed}{\end{equation}}
\newcommand{\bga}{\begin{eqnarray}}
\newcommand{\eda}{\end{eqnarray}}

\def\cbdu{\par{\raggedleft$\Box$\par}}

\newtheorem {Theorem}  {Theorem}

\numberwithin{Theorem}{section}

\newtheorem {Lemma}[Theorem]  {Lemma}

\theoremstyle{definition}

\theoremstyle{remark}
\newtheorem{Remark}[Theorem]{\bf Remark}

\expandafter\chardef\csname pre amssym.def
at\endcsname=\the\catcode`\@ \catcode`\@=11
\def\undefine#1{\let#1\undefined}
\def\newsymbol#1#2#3#4#5{\let\next@\relax
 \ifnum#2=\@ne\let\next@\msafam@\else
 \ifnum#2=\tw@\let\next@\msbfam@\fi\fi
 \mathchardef#1="#3\next@#4#5}
\def\mathhexbox@#1#2#3{\relax
 \ifmmode\mathpalette{}{\m@th\mathchar"#1#2#3}%
 \else\leavevmode\hbox{$\m@th\mathchar"#1#2#3$}\fi}
\def\hexnumber@#1{\ifcase#1 0\or 1\or 2\or 3\or 4\or 5\or 6\or 7\or 8\or
 9\or A\or B\or C\or D\or E\or F\fi}

\font\teneufm=eufm10 \font\seveneufm=eufm7 \font\fiveeufm=eufm5
\newfam\eufmfam
\textfont\eufmfam=\teneufm \scriptfont\eufmfam=\seveneufm
\scriptscriptfont\eufmfam=\fiveeufm

\catcode`\@=\csname pre amssym.def at\endcsname

\newcounter{remark}
\def  \12  {{\frac{1}{2}}}

\def\build#1_#2^#3{\mathrel{\mathop{\kern 0pt#1}\limits_{#2}^{#3}}}

\numberwithin{equation}{section}

\begin{document}
%\currannalsline{0}{2006}

%\title[Local existence for Hall-MHD]{Local existence for the non-resistive Hall-magneto-hydrodynamics system in $\R^n$}
\title[1D electron MHD]{Well-posedness and blowup of 1D electron magnetohydrodynamics}

%\author{hello}

\author [Mimi Dai]{Mimi Dai}

\address{Department of Mathematics, Statistics and Computer Science, University of Illinois at Chicago, Chicago, IL 60607, USA}
\email{mdai@uic.edu}

\thanks{M. Dai is partially supported by the NSF grant DMS--2308208 and Simons Foundation. }

\begin{abstract}

The one-dimensional toy models proposed for the three-dimensional electron magnetohydrodynamics in our previous work \cite{Dai-1d-emhd} share some similarities with the original dynamics under certain symmetry. We continue to study the well-posedness issue and explore the potential singularity formation scenario for these models.

%We propose one-dimensional reduced models for the three-dimensional electron magnetohydrodynamics which involves a highly nonlinear Hall term with intricate structure. The models contain nonlocal nonlinear terms which are more singular than that of the one-dimensional models for the Euler equation and the surface quasi-geostrophic equation. Local well-posedness is obtained in certain circumstances. Moreover, for a model with nonlocal transport term, we show that singularity develops in finite time for a class of initial data.

\bigskip


\hspace{0.02cm}CLASSIFICATION CODE: 35Q35, 76B03, 76D03, 76W05.
\end{abstract}

\maketitle

\section{Introduction}
%In order to capture the rapid magnetic reconnection phenomena in plasma physics, the authors of \cite{ADFL} rigorously derived the following incompressible magnetohydrodynamics (MHD) model with Hall effect 
%\begin{equation}\label{mhd}
%\begin{split}
%u_t+(u\cdot\nabla) u-(B\cdot\nabla) B+\nabla P=&\ \nu\Delta u, \\
%B_t+(u\cdot\nabla) B-(B\cdot\nabla) u +\nabla\times ((\nabla\times B)\times B)=&\ \mu\Delta B, \\
%\nabla\cdot u=0, \ \ \nabla\cdot B=&\ 0,
%\end{split}
%\end{equation}
%considered on the space time domain $\Omega\times [0,\infty)$ with $\Omega\subset \mathbb R^3$. In system (\ref{mhd}), the unknowns are the velocity field $u$, the magnetic field $B$ and the scalar pressure function $P$. The constant coefficients $\nu$ and $\mu$ denote respectively the kinetic viscosity and magnetic resistivity. The Hall term $\nabla\times ((\nabla\times B)\times B)$ appears to be more singular than other nonlinear terms in the system. Without it, system (\ref{mhd}) reduces to the classical MHD which shares many similarities with the Navier-Stokes equation (NSE). With the presence of the Hall term, (\ref{mhd}) no longer has a natural scaling as the classical MHD does. Instead, the Hall term introduces new scaling and geometry properties as well as new challenges. It seems desirable to understand the Hall term as a single target. 

The following electron magnetohydrodynamics (MHD) system 
\begin{equation}\label{emhd}
\begin{split}
B_t+ \nabla\times ((\nabla\times B)\times B)=&\ \mu\Delta B, \\
\nabla\cdot B=&\ 0
\end{split}
\end{equation}
arises in plasma physics phenomena where magnetic reconnection is extremely rapid and the slow ion flow motion can be neglected. In \eqref{emhd} the unknown $B: \mathbb R^3\times [0,\infty)\to \mathbb R^3$ denotes the magnetic field and the parameter $\mu\geq 0$ represents the resistivity. The great challenge of studying \eqref{emhd} comes from the nonlinear term $\nabla\times ((\nabla\times B)\times B)$, referred as the Hall term, which is rather singular. 
Since the mathematical derivation of \eqref{emhd} in \cite{ADFL}, a lot of efforts have been made to understand it as well as the full Hall MHD system from the analytical point of view, in terms of the basic questions of well-posedness, ill-posedness and possible finite time blowup, see \cite{CL, CS, CWW, CWeng, Dai22, Dai18, DO, JO, JO1, JO2} and references therein. 

In the literature we have seen extensive investigations of one-dimensional (1D) toy models for fluid equations such as the Euler/Navier-Stokes equations in the attempt to gain insights of their nonlinear structures, dating back to the early works \cite{CLM, DeG1, DeG2} and till the most recent studies \cite{Ch1, Ch2, Ch3, CHH, CHK, CKY, EGM, EJ, JSS, LLR, LH, OSW, Sak}. In particular, many articles have contributed to the topic of singularity formation for 1D simplified models of fluid equations, for instance see \cite{CCCF, CCF1, Do, DL, HLSWY, LR1, LR2} besides the ones mentioned above.
Inspired by the massive contributions, we proposed some nonlocal 1D models for the electron MHD in \cite{Dai-1d-emhd}. While the spatial geometric features of the original dynamics are suppressed in the 1D simplifications, the models do capture the two important nonlinear effects: transport and stretching. Well-posdness of the models with resistivity is obtained. Moreover, for the model with only transport effect and no resistivity, we discovered a scenario where a class of initial data yields solution that develops singularity formation in finite time.

To gain further insights of understanding the nonlinear structures encoded in the Hall term of the electron MHD, we continue to study the well-posedness issue and explore scenarios of finite time singularity formation for the 1D simplified models.
 
%\medskip

%\subsection{1D models for fluid equations} 
%As an attempt to understand the vorticity form of the Euler equation, a 1D model with nonlocal nonlinear structure was introduced by Constantin, Lax and Majda \cite{CLM}. This model was later generalized by De Gregorio \cite{DeG1, DeG2} and studied by many authors, for instance, see \cite{Ch1, Ch2, Ch3, CHH, CHK, CKY, EGM, EJ, LH, OSW, Sak}. Regarding these 1D toy models for the Euler equation and NSE, one important issue is the competition and balance between the stretching effect and the transport effect. With the efforts of the aforementioned authors, this question has been understood very well. Various well-posedess results were obtained in certain contexts, while solutions with finite time singularities were constructed in contrast scenarios. 

%It is known that the 2D SQG has strikingly similar features with the 3D Euler equation. In this vein of gaining insights from simplified models, 1D reduced models for the SQG were also proposed and investigated \cite{BLM, CasC, CCCF, CCF1, CCF2, Do, Dong, LR1, LR2, Mor, SV}. Both the problems of well-posedness and finite time blowup have been studied in depth. Readers with interest are referred to the works above. 

\medskip

\subsection{1D models for electron MHD} 
\label{sec-1dmhd}
% The toy models for the Hall MHD will be studied in a separate paper.
We recall the 1D nonlocal nonlinear models proposed for the electron MHD (\ref{emhd}) in \cite{Dai-1d-emhd}. 

Denoting $J=\nabla\times B$ by the current density, the Biot-Savart law 
\begin{equation}\notag%\label{BS}
B=\nabla\times (-\Delta)^{-1} J
\end{equation}
recovers $B$ from $J$. Since $\nabla \cdot B=0$, we can write the Hall term as 
\[\nabla\times ((\nabla\times B)\times B)=(B\cdot\nabla) J -(J\cdot\nabla) B.\]
Thus the electron MHD with generalized resistivity reads
\begin{equation}\label{emhd-1}
B_t+(B\cdot\nabla) J -(J\cdot\nabla) B+\mu \Lambda^\alpha B=0, \ \ \ \Lambda=(-\Delta)^{\frac12}
\end{equation}
with $\alpha\geq 0$.

In 1D setting, we apply the Hilbert transform
\begin{equation}\label{H1}
\mathcal H f=\frac{1}{\pi} P.V. \int_{-\infty}^{\infty}\frac{f(y)}{x-y}\, dy
\end{equation}
to define an analogy of the Biot-Savart law in the following way
\[B_x=\mathcal H J.\]
Given $J$, $B$ is uniquely determined under the zero-mean Gauge, i.e.
\begin{equation}\notag%\label{gauge1}
\int B(t,x)\,dx =0.
\end{equation}
The facts
 \[\Lambda f=\mathcal H\partial_x f=\partial_x\mathcal H f, \ \ \ \mathcal H\mathcal H f=-f\]
will be used often.  

The 1D toy model 
%on $[-\pi, \pi]$,
\begin{equation}\label{emhd-1d}
\begin{split}
B_t+B J_x-J B_x+\mu\Lambda^\alpha B=&\ 0,  \\
B_x=&\ \mathcal H J
\end{split}
\end{equation}
was introduced in \cite{Dai-1d-emhd} to mimic (\ref{emhd-1}).

We observe that if $B(x,t)$ is a solution to \eqref{emhd-1d} with initial data $B_0(x)$, the rescaled field $B_\lambda(x,t):=\lambda^{\alpha-2} B(\lambda x, \lambda^\alpha t)$ with a parameter $\lambda$ is a solution to \eqref{emhd-1d} with data $\lambda^{\alpha-2} B_0(\lambda x)$. The Sobolev space $\dot H^{\frac52-\alpha}$ is scaling invariant and referred critical for \eqref{emhd-1d}. % Note that $L^2$ is critical when $\alpha=\frac52$. By convention, (\ref{emhd-1d}) is referred to be supercritical, critical and subcritical respectively for $\alpha<\frac52$, $\alpha=\frac52$ and $\alpha>\frac52$. 

%In general belief, global well-posedness may be obtained for critical and subcritical systems, with $\alpha\geq \frac52$ in our context. However, the aforementioned lack of cancellation symmetry prevents us to achieve the global well-posedness by either using standard energy method or the approach of continuity of moduli. In another direction, since $BJ_x$ is more singular than $JB_x$, we conjecture that there exists initial data such that solutions of (\ref{e1d2}) with appropriate $\alpha$ blows up in finite time. These questions will be addressed in future work. 

\medskip

\subsection{Similarities between the 1D models and the axisymmetric electron MHD} 
\label{sec-compare}
We discuss in the following that the 1D model \eqref{emhd-1d} resembles the electron MHD in an axisymmetric setting. A point $(x,y,z)\in \mathbb R^3$ has the cylindrical coordinate $(r,\theta,z)$ with 
\[r=\sqrt{x^2+y^2}, \ \ \theta=\arctan (y/x).\]
The cylindrical coordinate basis $\{e_r, e_\theta, z\}$ is given by 
\[e_r=(\cos\theta, \sin\theta, 0), \ \ e_\theta=(-\sin\theta, \cos\theta, 0), \ \ e_z=(0,0,1).\]
We assume the magnetic field $B$ is axially symmetric, i.e. independent of the variable $\theta$, and in cylindrical coordinate it is written as
\[B(x,y,z,t)= b_r(r,z,t) e_r+b_\theta(r,z,t) e_\theta+b_z(r,z,t) e_z.\]
Hence we have for the current density 
\[J(x,y,z,t)=\nabla\times B(x,y,z,t)= j_r(r,z,t) e_r+j_\theta(r,z,t) e_\theta+j_z(r,z,t) e_z\]
with 
\begin{equation}\label{eq-j}
j_r=-\partial_z b_\theta, \ \ j_\theta=\partial_z b_r-\partial_r b_z, \ \ j_z=\frac1r\partial_r(rb_\theta).
\end{equation}
Converting to cylindrical coordinate, the axisymmetric electron MHD \eqref{emhd} with $\mu=0$ (non-resistive) is given by 
\begin{equation}\label{emhd-axis}
\begin{split}
\partial_t b_r-\partial_z(j_zb_r-j_rb_z)&=0,\\ %\mu \left(\frac1r \partial_r\left( r\partial_r b_r\right)+\partial_z^2 b_r-\frac{b_r}{r^2} \right),\\
\partial_t b_\theta+\partial_z(j_\theta b_z-j_zb_\theta)-\partial_r(j_rb_\theta-j_\theta b_r)&=0, \\ %\mu \left(\frac1r \partial_r\left( r\partial_r b_\theta\right)+\partial_z^2 b_\theta-\frac{b_\theta}{r^2} \right),\\
\partial_t b_z+\frac1r\partial_r\left(r(j_zb_r-j_rb_z)\right)&=0, \\ %\mu \left(\frac1r \partial_r\left( r\partial_r b_z\right)+\partial_z^2 b_z\right),\\
\partial_rb_r+\frac1rb_r+\partial_zb_z&=0.
\end{split}
\end{equation}

Below are some intriguing observations about the axisymmetric system \eqref{emhd-axis}. 
\begin{itemize}
\item [(i)] Vanishing swirl component: if the swirl component satisfies $b_\theta\equiv 0$, we have $j_r=j_z\equiv 0$ in view of \eqref{eq-j}.  It then follows from the first and third equations of \eqref{emhd-axis} that $\partial_t b_r=\partial_tb_z=0$; hence $b_r(r,z,t)=b_r(r,z,0)$ and $b_z(r,z,t)=b_z(r,z,0)$ for all $t\geq 0$. Although the dynamics is trivial in this case, some interesting stationary solutions may arise here. Furthermore, it indicates that the swirl component $b_\theta$ plays a dominant role in the dynamics. One may speculate that, if blowup occurs, it would have to be related to the growth of certain norms of $b_\theta$. 
\item [(ii)] Vanishing radial and vertical components: in contrast with (i), if the radial and vertical components vanish, i.e. $b_r\equiv 0$ and $b_z\equiv 0$, it follows from \eqref{eq-j} that $j_\theta\equiv 0$; consequently system \eqref{emhd-axis} reduces to one single equation of $b_\theta$,
\begin{equation}\notag
\partial_t b_\theta-\partial_z(j_zb_\theta)-\partial_r(j_rb_\theta)=0.
\end{equation}
Applying \eqref{eq-j} again, the equation above can be written as 
\begin{equation}\label{eq-b1}
\partial_t b_\theta-b_\theta(\partial_z j_z+\partial_r j_r)+\frac{b_\theta}{r}j_r=0.
\end{equation}
\item [(iii)] Resemblance of the swirl component and the 1D models: we note that equations \eqref{emhd-1d} and \eqref{eq-b1} share similarities among their nonlinear structures with $BJ_x$ corresponding to $b_\theta(\partial_z j_z+\partial_r j_r)$ and $B_xJ$ corresponding to $\frac{b_\theta}{r}j_r$.
If we further rewrite the ``stretching" term using \eqref{eq-j}
\[b_\theta(\partial_z j_z+\partial_r j_r)=-\frac{b_\theta}{r}j_r\]
which turns out to be the same as the transport term. It indicates that in the evolution of swirl component, the nonlinear structure has the only feature of transport. 
\end{itemize}

The observations above motivate our study of the 1D model \eqref{emhd-1d}. Although \eqref{emhd-1d} seems to be over simplified, it actually inherits the essential nonlinear structures of the original electron MHD.

\medskip

\subsection{Main results}
The local well-posedness of (\ref{emhd-1d}) with $\alpha>2$ in the Sobolev space $H^{s}(\mathbb R)$ with $s\geq \frac52-\alpha$ and smoothing estimates were obtained in \cite{Dai-1d-emhd}. The lack of divergence free property of $B$ in 1D setting is an obstacle to establish well-posedenss for weaker resistivity $\alpha\leq 2$. In this paper we further explore the nonlinear structures and discover mechanisms to move derivatives to low modes in the a prior estimates. We are able to improve the well-posedness to the case with resistivity $\alpha>1$. Since \eqref{emhd-1d} is subcritical for $\alpha>\frac52$, we only focus on the supercritical and critical case of $1<\alpha\leq \frac52$.
\begin{Theorem}\label{thm-local}
Let $\alpha\in(1, \frac52]$ and $B_0\in  H^{s}(\mathbb R)$ with $s\geq \frac52-\alpha$. There exists a time $T>0$ depending on $\|B_{0}\|_{H^s}$ and $\mu$ such that there exists a unique solution $B(x,t)$ to (\ref{emhd-1d}) with initial data $B(x,0)=B_0$ on $[0,T)$ satisfying
\begin{equation}\label{priori0}
B\in C\left([0,T);   H^{s}(\mathbb R)\right)\cap L^2\left([0,T);  H^{s+\frac{\alpha}2}(\mathbb R)\right)
\end{equation}
and 
\begin{equation}\label{est-b1}
\sup_{0<t<T} t^{\frac{s-\frac52+\alpha}{\alpha}} \|B(t)\|_{ \dot H^{s}(\mathbb R)}<\infty,
\end{equation}
\begin{equation}\label{est-b2}
\lim_{t\to 0} t^{\frac{s-\frac52+\alpha}{\alpha}} \|B(t)\|_{\dot H^{s}(\mathbb R)}=0, \ \ s>\frac52-\alpha.
\end{equation}
\end{Theorem}

\begin{Remark}\label{rek1}
We point out that
\begin{itemize}
\item [(i)]
To improve the case from $\alpha\in(2,\frac52]$ to $\alpha\in(1,\frac52]$, although we lose the divergence free property of $B$ and the associated cancelations in 1D setting, we explore the nonlinear mechanisms allowing us to move derivatives to low modes. See the estimate of $I_1$ in Section \ref{sec-local}.
\item [(ii)] We shall only show the well-posedness in $H^{s}$ for $s>\frac52-\alpha$ in Section \ref{sec-local}. With slightly more delicate analysis as in \cite{Dai-1d-emhd}, we can show local well-posedness in $H^{\frac52-\alpha}$. 
\end{itemize}
\end{Remark}

%{\color{blue}This may be improved to $\alpha\in(1, \frac52]$ by using a different integration by parts trick for the term $BJ_x$. }

%\medskip

%For $\alpha=1$ we show the local existence of analytic solution to \eqref{emhd-1d}. 
%\begin{Theorem}\label{thm-local1}
%Let $\alpha\in(1, \frac52]$ and $B_0\in  H^{\frac52-\alpha}(\mathbb S^1)$. The statements of Theorem \ref{thm-local} are true for model (\ref{e1d3}).
%\end{Theorem}

%For small $\alpha$, say $\alpha\leq 1$, we prove the local existence of analytic solutions to (\ref{e1d3}). When $\alpha=0$, we view $\mu=0$ automatically.

%\begin{Theorem}\label{thm-analytic}
%Let $\alpha\in[0, 1]$ and $\mu\geq 0$. Assume the initial data $B_0\in L^2(\mathbb R)$ is analytic. There exists a time $T>0$ such that equation (\ref{e1d3}) has a real analytic solution on $[0,T)$ with $B(x,0)=B_0(x)$.
%\end{Theorem}

%In the absence of $BJ_x$, with the help of maximum principle, we show finite time singularity occurs for (\ref{e1d3}) with $\mu=0$ and a class of initial data. 
%There might be some hope to show local well-posedness with the help of maximum principle. 

As reflected in the discussion of item (iii) in Subsection \ref{sec-compare}, it is meaningful to consider \eqref{emhd-1d} with the term stretching $BJ_x$ removed, i.e.
\begin{equation}\label{emhd-transport}
B_t-\Lambda B B_x+\mu\Lambda^\alpha B=0.
\end{equation}
For our preference, we replaced $B$ in \eqref{emhd-1d} by $-B$.
When $\mu=0$, it was shown in \cite{Dai-1d-emhd} that there exists a class of initial data such that the solution of \eqref{emhd-transport} develops singularity in finite time. Specifically, the initial profile provides a hyperbolic flow around the origin $x=0$ and $B_x(x,t)$ approaches infinity at $x=0$ as $t$ approaches the singular time. The singularity scenario relies on the maximum principle and the transport feature of \eqref{emhd-transport} with $\mu=0$.  For $\mu>0$, the maximum principle remains; however, the equation is no longer purely transport. Nevertheless, \eqref{emhd-transport} with $\alpha=1$ can be formulated into a transport equation for $\bar B(x,t)=B(x,t)-\mu x$. Benefited from this observation, we explore a different scenario of singularity formation for \eqref{emhd-transport} with $\mu>0$, where the quantity $\Lambda B_{x}$ blows up in finite time. Without loss of generality, we can take $\mu=1$.

\begin{Theorem}\label{thm-blow}
Assume the initial profile $B_0\in L^2(\mathbb R)\cap C^3(\mathbb R)$ decays fast to 0 for $|x|\gg 1$ and there exists a point $x_0\in \mathbb R$ such that 
\begin{equation}\label{initial}
\partial_x B_0(x_0)=1, \ \ \partial_{xx} B_0(x_0)=0, \ \ \Lambda\partial_{x}B_0(x_0)>0
\end{equation}
The solution $B(x,t)$ of (\ref{emhd-transport}) with $\mu=\alpha=1$ and the initial data $B_0$ develops singularity in the sense that $\Lambda B_{x}$ blows up at a finite time.
%Let $B(x,t)$ be a solution to (\ref{eq-s1}) with $\mu=0$ and the initial data $B_0$. Then the norm $\|B_x\|_{L^\infty}$ blows up in finite time.
\end{Theorem}

\begin{Remark}\label{rek2}
A few remarks about Theorem \ref{thm-blow} are in order.
\begin{itemize}
\item [(i)] We argue that the initial profile $B_0(x)$ satisfying the conditions of Theorem \ref{thm-blow} exists. For instance, if we first neglect the conditions of $B_0\in L^2(\mathbb R)$ and fast decay, we can choose $B_0(x)=\sin x$. Then
\[\partial_xB_0=\cos x, \ \ \partial_{xx}B_0=-\sin x\]
and 
\[\Lambda\partial_{x}B_0=\Lambda \cos x=\mathcal H\partial_x \cos x=-\mathcal H \sin x=\cos x\]
where the sign choice is consistent with the fact $\Lambda =\mathcal H\partial_x$. Taking $x_0=0$, we have
\[\partial_xB_0(0)=1, \ \ \partial_{xx}B_0(0)=0, \ \ \Lambda\partial_{x}B_0(0)=1>0.\]
Hence \eqref{initial} is satisfied. To fulfill the $L^2$ and fast decay requirements, we can modify the choice of the profile above as
\[B_0(x)=e^{-x^4}\sin x.\]
Then we have
\begin{equation}\notag
\begin{split}
\partial_xB_0=&\ e^{-x^4}\left(\cos x-4x^3\sin x\right),\\
\partial_{xx}B_0=&\ e^{-x^4}\left(-\sin x-12x^2\sin x+16x^6\sin x-8x^3\cos x\right).
\end{split}
\end{equation}
It is easy to check that 
\[\partial_xB_0(0)=1, \ \ \partial_{xx}B_0(0)=0\]
and $x_0=0$ is the maximum point of $\partial_xB_0$, i.e. $\partial_xB_0(x)<1$ for $x\neq 0$. We further verify that
\begin{equation}\notag
\begin{split}
\Lambda\partial_{x}B_0(0)=&\ \frac1\pi P.V. \int_{-\infty}^\infty\frac{\partial_{y}B_0(0)-\partial_yB_0(y)}{(0-y)^2} \, dy\\
=&\ \frac1\pi P.V. \int_{-\infty}^\infty\frac{1-\partial_yB_0(y)}{y^2} \, dy\\
>&\ 0.
\end{split}
\end{equation}
Therefore all the conditions on $B_0$ are satisfied.

More generically, we can choose $B_0\in L^2(\mathbb R)\cap C^3(\mathbb R)$ that decays fast to 0 for $|x|\gg 1$. Moreover, we require that there exists $x_0$ such that $x_0$ is the maximum point of $\partial_xB_0$ and $\partial_xB_0(x_0)=1$. Consequently, we know 
\[\partial_{xx}B_0(x_0)=0, \ \ \Lambda\partial_{x}B_0(x_0)>0. \]

\item [(ii)] The existence of an analytic solution under such conditions of the initial data was shown in \cite{Dai-1d-emhd}.
\item [(iii)] The regularity condition $C^3(\mathbb R)$ for the initial data $B_0$ to guarantee the existence of an analytic solution is likely not optimal. We do not pursue to identify the optimal regularity condition for the initial data in this paper.
\item [(iv)] In the non-resistive case of \eqref{emhd-transport} with $\alpha=0$ (essentially $\mu=0$), a scenario of singularity formation is found in \cite{Dai-1d-emhd}. In view of the singularity formation phenomena of \eqref{emhd-transport} with $\alpha=0$ and $\alpha=1$, it is reasonable to expect that blowup may occur for \eqref{emhd-transport} with $0<\alpha<1$. This will be addressed in future work.
%\item [(i)]
%One example of such initial profile functions satisfying the conditions of Theorem \ref{thm-s1} is 
%\[B_0(x)=e^{-x^2}\tan^{-1}x.\]
%\item [(ii)] The authors of \cite{SV} showed finite time singularity formation for the model 
%\begin{equation}\label{SV-model}
%\theta_t+ \theta_x (\Lambda^s\mathcal H\theta)=0, \ \ \ s\in(-1,1). 
%\end{equation}
%Note (\ref{e1d4}) is not a special case of (\ref{SV-model}); rather it is a ``boarderline'' situation of (\ref{SV-model}) with $s=1$ and an additional operator $\mathcal H$.
\end{itemize}
\end{Remark}

\bigskip

\section{Notations and preliminaries}

\subsection{Notations}
As usual, we use $C$ to denote a general constant which may be different at different lines. For simplicity, the notation $\lesssim$ is adapted for $\leq$ up to a constant.

\medskip

\subsection{Littlewood-Paley theory} 
To show the well-posedness, we apply the framework of Littlewood-Paley decomposition theory which is briefly recalled below. For an integer $q$ we
denote the frequency number $\lambda_q=2^q$. Take $\chi\in C_0^\infty(\mathbb R^n)$ to be a nonnegative function such that
\begin{equation}\notag
\chi(\xi)=
\begin{cases}
1, \ \ \ \mbox{for} \ \ |\xi|\leq \frac34\\
0, \ \ \ \mbox{for} \ \ |\xi|\geq 1.
\end{cases}
\end{equation}
Define $\varphi(\xi)=\chi(\frac{\xi}{2})-\chi(\xi)$ and 
\begin{equation}\notag
\varphi_q(\xi)=
\begin{cases}
\varphi(\lambda_q^{-1}\xi), \ \ \ \mbox{for} \ \ q\geq0\\
\chi(\xi), \ \ \ \ \ \ \ \  \mbox{for} \ \ q= -1. 
\end{cases}
\end{equation}
The $q$-th Littlewood-Paley projection of a tempered distribution vector field $u$ is defined as
%\[u_q(x)=\Delta_q u(x):= \sum_{k\in \mathbb Z^n} \widehat u(k)\varphi_q(k) e^{i\frac{2\pi}{L} k\cdot x}\]
\[\Delta_q u=\varphi_q*u.\]
The decomposition 
\[u=\sum_{q=-1}^\infty u_q\]
holds in the distributional sense. To ease notations, we write
\[\Delta_q u=u_q, \ \ \ \ u_{\leq Q}=\sum_{q\geq-1}^Q u_q, \ \ \ \ \widetilde {\Delta_q} u=u_{q-1}+u_q+u_{q+1}.\]
We point out that the Sobolev norm $\|u\|_{\dot H^s}$ is equivalent to 
\[\left(\sum_{q=-1}^\infty \lambda_q^{2s} \|u_q\|_{L^2}^2\right)^{\frac12}.\]

%\medskip

%\subsection{Auxiliary estimates}
%\label{sec-basic}
%We recollect some standard estimates regarding the operator $e^{-\mu t\Lambda^{\alpha}}$ in the following.
%\begin{Lemma}\label{le-m1} \cite{Miu}
%Let $\alpha>0$ and $f\in L^2$. There exist constants $C_1(\alpha), C_2(\alpha)>0$ such that 
%\begin{equation}\notag
%e^{-C_1(\alpha)\mu \lambda_q^{\alpha}t} \|f_q\|_{L^2}\leq \left\|e^{-\mu t\Lambda^{\alpha}}* f_q \right\|_{L^2}
%\leq e^{-C_2(\alpha)\mu \lambda_q^{\alpha}t} \|f_q\|_{L^2}.
%\end{equation}
%\end{Lemma}

%\begin{Lemma}\label{le-m2} \cite{Miu}
%Let $\alpha>0$, $s\geq 0$ and $f\in L^2$. There exists a constant $C(\alpha, s)>0$ such that 
%\begin{equation}\notag
%\begin{split}
%\sup_{t\in(0,\infty)}t^{\frac{s}{\alpha}} \left\|e^{-\mu t\Lambda^{\alpha}}* f \right\|_{H^s}\leq&\ C(\alpha, s) \|f\|_{L^2},\\
%\lim_{t\to 0}t^{\frac{s}{\alpha}} \left\|e^{-\mu t\Lambda^{\alpha}}* f \right\|_{H^s}=&\ 0, \ \ \mbox{for} \ \ s>0.
%\end{split}
%\end{equation}
%In addition if $s\in[0,\frac{\alpha}{2}]$, we have 
%\begin{equation}\notag
% \left\|e^{-\mu t\Lambda^{\alpha}}* f \right\|_{L^{\frac{\alpha}{s}}(0,T; H^s)}\leq C(\alpha, s, T,f) \|f\|_{L^2}
%\end{equation}
%where $C(\alpha, s, T,f)$ is a constant in the form 
%\[C(\alpha, s, T, f)=C_1(\alpha,s, f)\left(1-e^{-C_2(\alpha,s, f)\mu T}\right)^{\frac{s}{\alpha}}\]
%with another two constants $C_1(\alpha,s,f), C_2(\alpha,s,f)>0$ independent of the time $T$.

%\end{Lemma}

%\medskip

The commutator
\[[\Delta_q, f] g=\Delta_q(fg)-f g_q \]
satisfies the following estimate.
\begin{Lemma}\label{le-1} \cite{Dong}
%Let $f$ and $g$ be functions on $\mathbb T^n$. 
Let $n$ be the spatial dimension. Assume $r_0\geq 0$, $\frac{n}{2}\leq r_1<1+\frac{n}{2}$, $r_2<\frac{n}{2}$ and $r_0+r_1+r_2>0$. If $f\in \dot H^{r_1}\cap  \dot H^{r_0+r_1}$ and $g\in \dot H^{r_2}\cap  \dot H^{r_0+r_2}$, there exist a constant $C=C(n, r_0, r_1, r_2)>0$ and a sequence $\{c_q\}\in l^2$ with $\|c_q\|_{l^2}\leq 1$ such that 
\begin{equation}\notag
\left\|[\Delta_q, f] g\right\|_{\dot H^{r_0}}\leq C c_q\lambda_q^{-(r_1+r_2-\frac{n}{2})} \left(\|f\|_{\dot H^{r_0+r_1}}\|g\|_{\dot H^{r_2}}+\|f\|_{\dot H^{r_1}}\|g\|_{\dot H^{r_0+r_2}}\right).
\end{equation} 
% Let $m\geq 0$, $s_1<1+\frac{n}{2}$, $s_2<\frac{n}{2}$ and $m+s_1+s_2>0$. Assume $f\in  H^{s_1}\cap  H^{m+s_1}$ and $g\in  H^{s_2}\cap  H^{m+s_2}$. There exist a constant $C=C(m, n, s_1, s_2)>0$ and a sequence $\{c_q\}\in l^2$ with $\|c_q\|_{l^2}\leq 1$ such that for any $q\geq 0$  
%\begin{equation}\notag
%\left\|\widetilde {\Delta_q}[f, \Delta_q] g\right\|_{L^2}\leq C c_q\lambda_q^{-(m+s_1+s_2-\frac{n}{2})} \left(\|f\|_{H^{m+s_1}}\|g\|_{H^{s_2}}+\|f\|_{H^{s_1}}\|g\|_{H^{m+s_2}}\right).
%\end{equation} 
\end{Lemma}

We also have
\begin{Lemma}\label{le-comm} 
For any $1<r<\infty$,  there exists a constant $C$ depending on $r$ such that 
\begin{equation}\notag
\left\|[\Delta_q, f_{\leq p-2}] \nabla g_q\right\|_{L^r}\leq C \|\nabla f_{\leq p-2}\|_{L^\infty}\|g_q\|_{L^r}.
\end{equation} 
\end{Lemma}
This estimate was shown in \cite{Dai-hmhd-well} for $f$ being divergence free. Slight modification of the proof gives the estimate for general $f$ which is not necessarily divergence free.

%\medskip

%The next two lemmas are also import in the estimates later. 
%\begin{Lemma}\label{le-d2} \cite{Dong}
%Let $f$ and $g$ be functions on $\mathbb T^n$. 
%Let $m\geq 0$, $s_1<\frac{n}{2}$ and $s_2\in \mathbb R$. Assume $f\in H^{s_1}\cap H^{m+s_1}$ and $g\in H^{s_2}\cap  H^{m+s_2}$. There exist a constant $C=C(m, n, s_1, s_2)>0$ and a sequence $\{c_q\}\in l^2$ with $\|c_q\|_{l^2}\leq 1$ such that for any $q\geq 0$
%\begin{equation}\notag
%\left\|\widetilde {\Delta_q} (fg_q)\right\|_{L^2}\leq C c_q\lambda_q^{-(m+s_1+s_2-\frac{n}{2})} \left(\|f\|_{H^{m+s_1}}\|g\|_{H^{s_2}}+\|f\|_{H^{s_1}}\|g\|_{H^{m+s_2}}\right).
%\end{equation} 

%\end{Lemma}

%\medskip

%\begin{Lemma}\label{le-d3} \cite{RS}
%Let $s_1, s_2<\frac n2$ and $s_1+s_2>0$. If $f\in H^{s_1}$ and $g\in H^{s_2}$ we have
%\[\|fg\|_{H^{s_1+s_2-\frac{n}{2}}}\leq C \|f\|_{H^{s_1}} \|g\|_{H^{s_2}}\]
%for a constant $C=C(n, s_1, s_2)>0$.
%\end{Lemma}

\medskip

%\subsection{Properties of Hilbert transform}
%The Hilbert transform has the following simple properties
%\begin{equation}\notag
%\begin{split}
%\mathcal H(cf)=&\ c\mathcal Hf, \ \ \mbox{for a constant} \ \ c,\\
%\mathcal H\sin(kx)=&-\cos(kx),\ \ \ \mathcal H\cos(kx)=\sin(kx).
%\end{split}
%\end{equation}
%And more generally, we have
%\begin{equation}\notag
%\mathcal H\sin(kx+\theta)=-\cos(kx+\theta), \ \ \mathcal H\cos(kx+\theta)=\sin(kx+\theta).
%\end{equation}

%\begin{equation}\notag
%\begin{split}
%\mathcal H(\mathcal H f)=&-f,\\
%\mathcal H(fg)=&\ f\mathcal H g+g\mathcal H f+\mathcal H(\mathcal Hf\mathcal Hg),\\
%\mathcal H(f\mathcal H f)=&\ \frac12\left((\mathcal H f)^2-f^2\right), \\
%\mathcal H(e^{ikx})=&\ i \cdot\mathrm{sign}(k) e^{ikx}.
%\end{split}
%\end{equation}

%For any periodic function $f$, the mean value of its Hilbert transform is zero, that is
%\begin{equation}\label{Hm}
%\int \mathcal Hf\, dx=0.
%\end{equation}

Denote the Riesz potential operator $\mathcal I_r$ by
\begin{equation}\notag
\mathcal I_r f(x)=\int_{\mathbb R^n} \frac{f(y)}{|x-y|^{n-r}}\, dy, \ \ \ 0<r<n.
\end{equation} 
We observe that the Riesz potential $\mathcal I_r$ is inverse to the fractional Laplacian operator $\Lambda ^{r}$.
We have the following sharp commutator estimate of the Coifman-Meyer type.

\begin{Lemma}\label{le-2}\cite{LS} 
Let $\gamma\in (0,1]$ and $1<r, r_1, r_2<\infty$ with $\frac1{r_1}+\frac1{r_2}=\frac1r$. For $\sigma\in[\gamma,1)$, there exists a constant $C>0$ depending on the parameters $\gamma, r, r_1,r_2$ and $\sigma$ such that
\begin{equation}\notag
\|[\Lambda^\gamma, f]g\|_{L^r}\leq C \|\Lambda^{\sigma} f\|_{L^{r_1}}\|\mathcal I_{\sigma-\gamma} g\|_{L^{r_2}}.
\end{equation}
\end{Lemma}

In the end we recall that the Hilbert transform is a bounded linear operator from space $L^p$ to $L^p$. 

\begin{Lemma}\label{le-3} \cite{Zy} 
Let $1<p<\infty$. There exists a constant $C>0$ depending on $p$ such that
\begin{equation}\notag%\label{HLp}
\|\mathcal Hf\|_{L^p}\leq C(p) \|f\|_{L^p}.
\end{equation}
\end{Lemma}

\bigskip

\section{Local well-posedness}
\label{sec-local}

We show the well-posedness result of Theorem \ref{thm-local} for $1<\alpha\leq \frac52$ in this section. As customary, the main step is to establish the a priori estimate in the space $H^s$ for $s>\frac52-\alpha$; the existence then follows from a standard procedure of producing a sequence of approximating solutions and taking the limit of the sequence (or a subsequence). 
%This is achieved in three steps: we first establish the a priori estimates; we then solve the approximating system, see (\ref{sys-app}), to obtain a sequence of approximating solutions satisfying the a priori estimates; in the end, we show the sequence is a Cauchy sequence in appropriate functional spaces by employing the a priori estimates and hence converges to a solution. The main difficulties arise in the first step. 
In the 1D setting, the divergence free property $\nabla\cdot B=0$ is not valid. It causes several obstacles in our analysis. One immediate drawback is that we can not obtain the estimate of $B$ in $L^2$ in a trivial way as for the original 3D electron MHD system. We adapt the backward bootstrap argument from \cite{Dai-1d-emhd, Dong} to first establish an estimate in $\dot H^s$ for $s>\frac52-\alpha$ and then the estimate of $L^2$. The more essential difficulty is the lack of cancellations due to the absence of divergence free. To overcome it, as explained in Remark \ref{rek1} (i), we explore the mechanisms of moving derivatives from high modes to low ones. More details are demonstrated in the estimate of $I_1$ below.

%\medskip

\subsection{A priori estimate in $\dot H^s$ for $s>\frac52-\alpha$}
\label{sec-priori}
Acting the projection $\Delta_q$ on equation (\ref{emhd-1d}), multiplying the resulted equation by $B_q$ and integrating yields
\begin{equation}\notag
\begin{split}
\frac12\frac{d}{dt}\|B_q\|_{L^2}^2+\mu \|\Lambda^{\frac\alpha2}B_q\|_{L^2}^2
=&-\int_{\mathbb R}(BJ_x)_qB_q dx+\int_{\mathbb R}(JB_x)_qB_q dx\\
=&\int_{\mathbb R}(BJ)_qB_{x,q} dx+2\int_{\mathbb R}(JB_x)_qB_q dx\\
=&-\int_{\mathbb R}(B\Lambda B)_qB_{x,q} dx-2\int_{\mathbb R}(\Lambda BB_x)_qB_q dx\\
\end{split}
\end{equation}
where we used integration by parts and the fact $J=-\Lambda B$. Multiplying the equation by $\lambda_q^{2s}$ and taking summation over $q\geq -1$ yields
\begin{equation}\label{est-energy}
\begin{split}
&\frac12\frac{d}{dt}\sum_{q\geq -1}\lambda_q^{2s}\|B_q\|_{L^2}^2+\mu\sum_{q\geq -1}\lambda_q^{2s} \|\Lambda^{\frac\alpha2}B_q\|_{L^2}^2\\
=&-\sum_{q\geq -1}\lambda_q^{2s}\int_{\mathbb R}(B\Lambda B)_qB_{x,q} dx-2\sum_{q\geq -1}\lambda_q^{2s}\int_{\mathbb R}(\Lambda BB_x)_qB_q dx\\
=&: -I-2 K.
\end{split}
\end{equation}
%\begin{equation}\notag
%\partial_t B_q+\Delta_q(BJ_x)-\Delta_q(JB_x)+\mu\Lambda^\alpha B_q=0.
%\end{equation}
%Applying the commutation notation we have
%\begin{equation}\notag
%\partial_t B_q+BJ_{x,q}-JB_{x,q}+\mu\Lambda^\alpha B_q=[B, \Delta_q]J_x+[J, \Delta_q]B_x.
%\end{equation}
%Multiplying the equation above by $B_q$ and integrating over $\mathbb S^1$ gives
%\begin{equation}\notag
%\begin{split}
%&\frac12\frac{d}{dt}\|B_q\|_{L^2}^2+\mu\lambda_q^\alpha \|B_q\|_{L^2}^2\\
%\leq &\ a\int_{\mathbb S^1} [B, \Delta_q]J_x B_q\, dx+b\int_{\mathbb S^1} [J, \Delta_q]B_x B_q\, dx\\
%&-a\int_{\mathbb S^1} BJ_{x,q} B_q\, dx-b\int_{\mathbb S^1} JB_{x,q} B_q\, dx\\
%=&\ a \int_{\mathbb S^1} \widetilde{\Delta_q}\left([B, \Delta_q]J_x\right) B_q\, dx+b\int_{\mathbb S^1} \widetilde{\Delta_q}\left([J, \Delta_q]B_x\right) B_q\, dx\\
%&-a\int_{\mathbb S^1} \widetilde{\Delta_q}\left(BJ_{x,q}\right) B_q\, dx-b\int_{\mathbb S^1} \widetilde{\Delta_q}\left(JB_{x,q}\right) B_q\, dx.
%\end{split}
%\end{equation}

\subsection{Estimate of $I$}
The integral $I$ can be decomposed using Bony's paraproduct as
\begin{equation}\notag
\begin{split}
I=&\sum_{q\geq -1}\sum_{|p-q|\leq 2}\lambda_q^{2s}\int_{\mathbb R}\Delta_q(B_{\leq p-2}\Lambda B_{p})B_{x,q} dx\\
&+\sum_{q\geq -1}\sum_{|p-q|\leq 2}\lambda_q^{2s}\int_{\mathbb R}\Delta_q(B_{p}\Lambda B_{\leq p-2})B_{x,q} dx\\
&+\sum_{q\geq -1}\sum_{p\geq q-2}\lambda_q^{2s}\int_{\mathbb R}\Delta_q(B_{p}\Lambda \widetilde B_{p})B_{x,q} dx\\
=&: I_{1}+I_{2}+I_{3}.
\end{split}
\end{equation}
We would like to move some derivative to the low modes part $B_{\leq p-2}$ in $I_1$. In order to do so, we apply a commutator to rearrange $I_1$ as
\begin{equation}\notag
\begin{split}
I_1=&\sum_{q\geq -1}\sum_{|p-q|\leq 2}\lambda_q^{2s}\int_{\mathbb R}[\Delta_q,B_{\leq p-2}]\Lambda B_{p}B_{x,q} dx\\
&+\sum_{q\geq -1}\sum_{|p-q|\leq 2}\lambda_q^{2s}\int_{\mathbb R}B_{\leq q-2}\Delta_q\Lambda B_{p}B_{x,q} dx\\
&+\sum_{q\geq -1}\sum_{|p-q|\leq 2}\lambda_q^{2s}\int_{\mathbb R}\left((B_{\leq p-2}-B_{\leq q-2})\Delta_q\Lambda B_{p}\right)B_{x,q} dx\\
=&:I_{11}+I_{12}+I_{13}.
\end{split}
\end{equation}
Applying the fact $\Lambda=\mathcal H\partial_x=\partial_x\mathcal H$, H\"older's inequality, the commutator estimate in Lemma \ref{le-comm}, the boundedness of Hilbert transform in Lemma \ref{le-3} and Bernstein's inequality, we infer
\begin{equation}\label{est-I11}
\begin{split}
|I_{11}|=&\left|\sum_{q\geq -1}\sum_{|p-q|\leq 2}\lambda_q^{2s}\int_{\mathbb R}[\Delta_q,B_{\leq p-2}]\partial_x\mathcal HB_{p}B_{x,q} dx\right|\\
\leq&\sum_{q\geq -1}\sum_{|p-q|\leq 2}\lambda_q^{2s}\|[\Delta_q,B_{\leq p-2}]\partial_x\mathcal HB_{p}\|_{L^2}\|B_{x,q}\|_{L^2}\\
\lesssim&\sum_{q\geq -1}\sum_{|p-q|\leq 2}\lambda_q^{2s}\|\partial_xB_{\leq p-2}\|_{L^\infty}\|\mathcal HB_{p}\|_{L^2}\|B_{x,q}\|_{L^2}\\
\lesssim&\sum_{q\geq -1}\lambda_q^{2s+1}\|B_{q}\|_{L^2}^2\sum_{p\leq q}\|\partial_xB_{p}\|_{L^\infty}\\
\lesssim&\sum_{q\geq -1}\lambda_q^{2s+1}\|B_{q}\|_{L^2}^2\sum_{p\leq q}\lambda_p^{\frac32}\|B_{p}\|_{L^2}.
\end{split}
\end{equation}
Next we apply Young's inequality and Jensen's inequality for some $0<\delta_1<2$ to obtain
\begin{equation}\label{est-I11-b}
\begin{split}
|I_{11}|
\lesssim & \sum_{q\geq -1}\left(\lambda_q^{(s+\frac{\alpha}2)\delta_1}\|B_q\|_{L^2}^{\delta_1}\right)\left(\lambda_q^{s(2-\delta_1)}\|B_q\|_{L^2}^{2-\delta_1}\right) \sum_{p\leq q}\lambda_p^{s}\|B_{p}\|_{L^2}  \lambda_{p-q}^{\frac32-s}\lambda_q^{\frac52-\frac{\alpha}2\delta_1-s}\\
\leq&\ \frac{1}{64}\mu\sum_{q\geq -1}\lambda_q^{2s+\alpha}\|B_q\|_{L^2}^2\\
&+\frac{C}{\mu}\sum_{q\geq -1}\lambda_q^{2s}\|B_q\|_{L^2}^2\left(\sum_{p\leq q}\lambda_p^{s}\|B_{p}\|_{L^2}  \lambda_{p-q}^{\frac32-s}\lambda_q^{\frac52-\frac{\alpha}2\delta_1-s} \right)^{\frac2{2-\delta_1}}\\
\leq&\ \frac{1}{64}\mu\sum_{q\geq -1}\lambda_q^{2s+\alpha}\|B_q\|_{L^2}^2+\frac{C}{\mu}\left(\sum_{q\geq -1}\lambda_q^{2s}\|B_q\|_{L^2}^2\right)^{1+\frac1{2-\delta_1}}
\end{split}
\end{equation}
provided the parameters satisfy
\begin{equation}\label{para1}
\frac{3}{2}-s>0, \ \ \ \frac52-\frac{\alpha}2\delta_1-s\leq 0.
\end{equation}
Note that, if $\alpha>1$, there is $\delta_1\in(0,2)$ such that \eqref{para1} is satisfied for $\frac52-\alpha< s<\frac32$. For instance, one can choose $\delta_1\in[\frac{5-2s}{\alpha},2)$. 

Invoking the fact $\sum_{|p-q|\leq 2} \Delta_q\Lambda B_p=\Delta_q \Lambda B=\Lambda B_q$, the term $I_{12}$ can be further rewritten as
\begin{equation}\notag
\begin{split}
I_{12}=&\sum_{q\geq -1}\lambda_q^{2s}\int_{\mathbb R}B_{\leq q-2}\Lambda B_{q}B_{x,q} dx\\
=&\sum_{q\geq -1}\lambda_q^{2s}\int_{\mathbb R}\Lambda^{\frac12}(B_{\leq q-2}B_{x,q})\Lambda^{\frac12} B_q dx\\
=&\sum_{q\geq -1}\lambda_q^{2s}\int_{\mathbb R}[\Lambda^{\frac12}, B_{\leq q-2}]B_{x,q}\Lambda^{\frac12} B_q dx\\
&+\sum_{q\geq -1}\lambda_q^{2s}\int_{\mathbb R} B_{\leq q-2}\Lambda^{\frac12}  B_{x,q}\Lambda^{\frac12} B_q dx
\end{split}
\end{equation}
where we used integration by parts. One more time of integration by parts gives us
\begin{equation}\notag
\begin{split}
I_{12}=&\sum_{q\geq -1}\lambda_q^{2s}\int_{\mathbb R}[\Lambda^{\frac12}, B_{\leq q-2}]B_{x,q}\Lambda^{\frac12} B_q dx\\
&+\frac12\sum_{q\geq -1}\lambda_q^{2s}\int_{\mathbb R} B_{\leq q-2}\partial_x\left(\Lambda^{\frac12}  B_{q}\Lambda^{\frac12} B_q\right) dx\\
=&\sum_{q\geq -1}\lambda_q^{2s}\int_{\mathbb R}[\Lambda^{\frac12}, B_{\leq q-2}]B_{x,q}\Lambda^{\frac12} B_q dx\\
&-\frac12\sum_{q\geq -1}\lambda_q^{2s}\int_{\mathbb R} B_{x,\leq q-2} \Lambda^{\frac12} B_q\Lambda^{\frac12}  B_{q}dx\\
=&:I_{121}+I_{122}.
\end{split}
\end{equation}
Applying H\"older's inequality and the commutator estimate in Lemma \ref{le-2} with $r=2$, $\gamma=\frac12$, $\frac12\leq \sigma<1$ and $1<r_1,r_2<\infty$ satisfying $\frac1{r_1}+\frac1{r_2}=\frac12$, we deduce
\begin{equation}\notag
\begin{split}
|I_{121}|\leq & \sum_{q\geq -1}\lambda_q^{2s}\|[\Lambda^{\frac12}, B_{\leq q-2}]B_{x,q}\|_{L^2}\|\Lambda^{\frac12}B_q\|_{L^2}\\
\lesssim & \sum_{q\geq -1}\lambda_q^{2s}\|\Lambda^{\sigma} B_{\leq q-2}\|_{L^{r_1}}\|\mathcal I_{\sigma-\frac12}B_{x,q}\|_{L^{r_2}}\|\Lambda^{\frac12}B_q\|_{L^2}\\
\lesssim & \sum_{q\geq -1}\lambda_q^{2s}\|\Lambda^{\sigma} B_{\leq q-2}\|_{L^{r_1}}\|\Lambda^{\frac12-\sigma}B_{x,q}\|_{L^{r_2}}\|\Lambda^{\frac12}B_q\|_{L^2}.
\end{split}
\end{equation}
Choosing $\sigma=1-\varepsilon$ and $r_2=2+\varepsilon$ for arbitrarily small $\varepsilon>0$ and using Bernstein's inequality gives
\begin{equation}\notag
\begin{split}
|I_{121}|
\lesssim & \sum_{q\geq -1}\lambda_q^{2s}\|\Lambda^{\sigma} B_{\leq q-2}\|_{L^{r_1}}\|\Lambda^{\frac12-\sigma}B_{x,q}\|_{L^{r_2}}\|\Lambda^{\frac12}B_q\|_{L^2}\\
\lesssim & \sum_{q\geq -1}\lambda_q^{2s+2-\sigma+\frac{\varepsilon}{2(2+\varepsilon)}}\|B_q\|_{L^2}^2\sum_{p\leq q-2}\lambda_p^{\sigma+\frac1{2+\varepsilon}}\|B_{p}\|_{L^{2}}\\
\lesssim & \sum_{q\geq -1}\lambda_q^{2s+1+\varepsilon+\frac{\varepsilon}{2(2+\varepsilon)}}\|B_q\|_{L^2}^2\sum_{p\leq q-2}\lambda_p^{1-\varepsilon+\frac1{2+\varepsilon}}\|B_{p}\|_{L^{2}}.
\end{split}
\end{equation}
Then we apply Young's inequality and Jensen's inequality for some $\delta_2\in(0,2)$ to obtain
\begin{equation}\notag
\begin{split}
|I_{121}|
\lesssim & \sum_{q\geq -1}\left(\lambda_q^{(s+\frac{\alpha}2)\delta_2}\|B_q\|_{L^2}^{\delta_2}\right)\left(\lambda_q^{s(2-\delta_2)}\|B_q\|_{L^2}^{2-\delta_2}\right) \sum_{p\leq q-2}\lambda_p^{s}\|B_{p}\|_{L^2} \\
&\cdot  \lambda_{p-q}^{1-s-\varepsilon+\frac1{2+\varepsilon}} \lambda_q^{s-\frac{\alpha}2\delta_2-\frac1{2+\varepsilon}+\frac{\varepsilon}{2(2+\varepsilon)}+2\varepsilon}\\
\leq&\ \frac{1}{64}\mu\sum_{q\geq -1}\lambda_q^{2s+\alpha}\|B_q\|_{L^2}^2\\
&+\frac{C}{\mu}\sum_{q\geq -1}\lambda_q^{2s}\|B_q\|_{L^2}^2\left(\sum_{p\leq q-2}\lambda_p^{s}\|B_{p}\|_{L^2} \lambda_{p-q}^{1-s-\varepsilon+\frac1{2+\varepsilon}} \lambda_q^{s-\frac{\alpha}2\delta_2-\frac1{2+\varepsilon}+\frac{\varepsilon}{2(2+\varepsilon)}+2\varepsilon} \right)^{\frac2{2-\delta_2}}\\
\leq&\ \frac{1}{64}\mu\sum_{q\geq -1}\lambda_q^{2s+\alpha}\|B_q\|_{L^2}^2+\frac{C}{\mu}\left(\sum_{q\geq -1}\lambda_q^{2s}\|B_q\|_{L^2}^2\right)^{1+\frac1{2-\delta_2}}
\end{split}
\end{equation}
provided the parameters satisfy
\begin{equation}\label{para-2}
\begin{split}
1-s-\varepsilon+\frac1{2+\varepsilon}>0,\\
s-\frac{\alpha}2\delta_2-\frac1{2+\varepsilon}+\frac{\varepsilon}{2(2+\varepsilon)}+2\varepsilon\leq 0.
\end{split}
\end{equation}
Since $s>\frac52-\alpha$ and $\alpha>1$, we can choose $\varepsilon>0$ small enough such that 
\[\frac32-\alpha< \frac1{2+\varepsilon}-\varepsilon, \ \ \frac52-2\alpha<\frac1{2+\varepsilon}-\frac{\varepsilon}{2(2+\varepsilon)}-2\varepsilon, \]
and then the conditions of \eqref{para-2} will be satisfied for some $0<\delta_2<2$. 

With one derivative on the low modes in $I_{122}$, the estimate of $I_{122}$ is much easier. We apply H\"older's inequality and Bernstein's inequality to obtain
\begin{equation}\notag
\begin{split}
|I_{122}|\leq& \sum_{q\geq-1}\frac12\lambda_q^{2s}\|\Lambda^{\frac12}B_q\|_{L^2}^2\|B_{x,\leq q-2}\|_{L^\infty}\\
\lesssim& \sum_{q\geq-1}\lambda_q^{2s+1}\|B_q\|_{L^2}^2\sum_{p\leq q-2}\lambda_p^{\frac32}\|B_{p}\|_{L^2}\\
\end{split}
\end{equation}
which is in a similar form as that of $I_{11}$ in \eqref{est-I11}. Thus $I_{122}$ shares the same estimate of $I_{11}$ as in \eqref{est-I11-b}
\begin{equation}\notag
|I_{122}|
\leq \frac{1}{64}\mu\sum_{q\geq -1}\lambda_q^{2s+\alpha}\|B_q\|_{L^2}^2+\frac{C}{\mu}\left(\sum_{q\geq -1}\lambda_q^{2s}\|B_q\|_{L^2}^2\right)^{1+\frac1{2-\delta_1}}.
\end{equation}

Observing that the sum $\sum_{|p-q|\leq 2}(B_{\leq p-2}-B_{\leq q-2})$ contains only a few non-zero modes, the integral $I_{13}$ involves with only one summation of infinitely many terms. The estimate of $I_{13}$ is given by
\begin{equation}\label{est-I13}
\begin{split}
|I_{13}|\leq &\sum_{q\geq -1}\sum_{|p-q|\leq 2}\lambda_q^{2s}\|B_{\leq p-2}-B_{\leq q-2}\|_{L^\infty}\|\Delta_q\Lambda B_p\|_{L^2}\| B_{x,q}\|_{L^2}\\
\lesssim &\sum_{q\geq -1}\lambda_q^{2s+2}\|B_q\|_{L^\infty}\| B_{q}\|_{L^2}^2\\
\lesssim &\sum_{q\geq -1}\lambda_q^{2s+\frac52}\| B_{q}\|_{L^2}^3.
\end{split}
\end{equation}
We continue the estimate by using Young's inequality for some $0<\delta_3<2$
\begin{equation}\label{est-I13-b}
\begin{split}
|I_{13}|
\lesssim &\sum_{q\geq -1}\left(\lambda_q^{(s+\frac\alpha2)\delta_3}\|B_q\|_{L^2}^{\delta_3}\right)\left(\lambda_q^{s(3-\delta_3)}\|B_q\|_{L^2}^{3-\delta_3}\right)\lambda_q^{\frac52-s-\frac\alpha2\delta_3}\\
\leq &\ \frac{1}{64}\mu\sum_{q\geq -1}\lambda_q^{2s+\alpha}\|B_q\|_{L^2}^2+\frac{C}{\mu} \sum_{q\geq -1}\lambda_q^{\frac{2s(3-\delta_3)}{2-\delta_3}}\|B_q\|_{L^2}^{\frac{2(3-\delta_3)}{2-\delta_3}}\\
\leq &\ \frac{1}{64}\mu\sum_{q\geq -1}\lambda_q^{2s+\alpha}\|B_q\|_{L^2}^2+\frac{C}{\mu}\left( \sum_{q\geq -1}\lambda_q^{2s}\|B_q\|_{L^2}^2\right)^{\frac{3-\delta_3}{2-\delta_3}}\\
\end{split}
\end{equation}
provided that $\frac52-s-\frac\alpha2\delta_3\leq 0$. Note for $s>\frac52-\alpha$ and $0<\alpha\leq \frac52$, we can choose $\delta_3\in[\frac{5-2s}{\alpha}, 2)$.

To estimate $I_{2}$, 
we first apply H\"older's and Bernstein's inequalities
\begin{equation}\notag
\begin{split}
|I_{2}|
\leq &\sum_{q\geq -1}\sum_{|p-q|\leq 2}\lambda_q^{2s}\|B_p\|_{L^2}\|\Lambda B_{\leq p-2}\|_{L^\infty}\|B_{x,q}\|_{L^2}   \\
\lesssim &\sum_{q\geq -1}\lambda_q^{2s+1}\|B_q\|_{L^2}^2\sum_{p\leq q}\lambda_p^{\frac32}\|B_{p}\|_{L^2}   \\
\end{split}
\end{equation}
which is also similar to $I_{11}$. Hence we have
\begin{equation}\notag
|I_{2}|
\leq \frac{1}{64}\mu\sum_{q\geq -1}\lambda_q^{2s+\alpha}\|B_q\|_{L^2}^2+\frac{C}{\mu}\left(\sum_{q\geq -1}\lambda_q^{2s}\|B_q\|_{L^2}^2\right)^{1+\frac1{2-\delta_1}}.
\end{equation}

The estimate of $I_{3}$ is more involved. We start by rewriting it using a commutator
\begin{equation}\notag
\begin{split}
I_{3}=& \sum_{q\geq -1}\sum_{p\geq q-2}\lambda_q^{2s}\int_{\mathbb R}B_{p}\Delta_q \Lambda \widetilde B_{p} B_{x,q} dx\\
&+\sum_{q\geq -1}\sum_{p\geq q-2}\lambda_q^{2s}\int_{\mathbb R}[\Delta_q, B_{p}] \Lambda \widetilde B_{p}B_{x,q} dx\\
=&: I_{31}+I_{32}.
\end{split}
\end{equation}
Observing that $\Delta_q\Lambda \widetilde B_{p}=0$ for $|p-q|\geq 5$, the summation $\sum_{p\geq q-2}$ in $I_{31}$ has essentially a few non-zero terms. Thus
\begin{equation}\notag
\begin{split}
|I_{31}|\lesssim &\sum_{q\geq -1}\sum_{p\sim q}\lambda_q^{2s}\|\Lambda \widetilde B_p\|_{L^2}\|B_{p}\|_{L^2}\|B_{x,q}\|_{L^\infty}\\
\lesssim &\sum_{q\geq -1}\lambda_q^{2s+\frac52}\|B_q\|_{L^2}^3\\
\end{split}
\end{equation}
which is similar to $I_{13}$ as shown in \eqref{est-I13}. Therefore we have the similar estimate as in \eqref{est-I13-b} for $I_{31}$,
\begin{equation}\notag
|I_{31}|\leq  \frac{1}{64}\mu\sum_{q\geq -1}\lambda_q^{2s+\alpha}\|B_q\|_{L^2}^2+\frac{C}{\mu}\left( \sum_{q\geq -1}\lambda_q^{2s}\|B_q\|_{L^2}^2\right)^{\frac{3-\delta_3}{2-\delta_3}}.
\end{equation}

To estimate $I_{32}$, we apply H\"older's inequality and Lemma \ref{le-1} with $r_0=0$, $n=1$, $\frac12\leq r_1<\frac32$, $r_2<\frac12$ and $r_1+r_2>0$
\begin{equation}\notag
\begin{split}
|I_{32}|\lesssim & \sum_{q\geq -1}\lambda_q^{2s}\|B_{x,q}\|_{L^2}\sum_{p\geq q-2}\|[\Delta_q, B_{p}]\Lambda \widetilde B_{p}\|_{L^2}\\
\lesssim & \sum_{q\geq -1}c_q\lambda_q^{2s+1-(r_1+r_2-\frac12)}\|B_q\|_{L^2}\sum_{p\geq q-2}\|B_{p}\|_{\dot H^{r_1}}\|\Lambda \widetilde B_{p}\|_{\dot H^{r_2}}\\
\lesssim & \sum_{p\geq -1}\|B_{p}\|_{\dot H^{r_1}}\|\Lambda \widetilde B_{p}\|_{\dot H^{r_2}}\sum_{q\leq p+2}c_q\lambda_q^{2s+1-(r_1+r_2-\frac12)}\|B_q\|_{L^2}\\
\lesssim & \sum_{p\geq -1}\lambda_p^{1+r_1+r_2}\| B_{p}\|_{L^2}^2\sum_{q\leq p+2}c_q\lambda_q^{2s+\frac32-(r_1+r_2)}\|B_q\|_{L^2}.
\end{split}
\end{equation}
We continue to rearrange the terms in the last line for some $\delta_4\in(0,2)$ as 
\begin{equation}\notag
\begin{split}
|I_{32}|
%\lesssim & \sum_{p\geq -1}\lambda_p^{1+r_1+r_2}\| B_{p}\|_{L^2}^2\sum_{q\leq p+2}c_q\lambda_q^{2s+\frac32-(r_1+r_2)}\|B_q\|_{L^2}\\
\lesssim & \sum_{p\geq -1}\left(\lambda_p^{(s+\frac\alpha2)\delta_4}\| B_{p}\|_{L^2}^{\delta_4}\right)\left(\lambda_p^{s(2-\delta_4)}\| B_{p}\|_{L^2}^{2-\delta_4}\right)  \sum_{q\leq p+2}c_q\lambda_q^{s}\|B_q\|_{L^2}\\
&\cdot \lambda_q^{\frac52-s-\frac\alpha2\delta_4}\lambda_{p-q}^{1+r_1+r_2-2s-\frac\alpha2\delta_4}.
\end{split}
\end{equation}
We need to choose parameters satisfying 
\begin{equation}\label{para-3}
\begin{split}
\frac52-s-\frac\alpha2\delta_4&\leq 0,\\
1+r_1+r_2-2s-\frac\alpha2\delta_4&<0
\end{split}
\end{equation}
such that we can apply Young's and Jensen's inequalities.  For $0<\alpha\leq \frac52$ and $s>\frac52-\alpha$, we can choose $\delta_4\in[\frac{5-2s}{\alpha}, 2)$ and then the first condition of \eqref{para-3} is satisfied. We choose $r_1=\frac12$, $r_2=-\frac12+\varepsilon$ for a small constant $\varepsilon>0$ and $\delta_4$ close enough to 2. Then the second condition of \eqref{para-3} is also satisfied for $0<\alpha\leq \frac52$ and $s>\frac52-\alpha$.
With such choice of parameters we proceed the estimate by applying Young's inequality and Jensen's inequality 
\begin{equation}\notag
\begin{split}
|I_{32}|
\lesssim & \sum_{p\geq -1}\left(\lambda_p^{(s+\frac\alpha2)\delta_4}\| B_{p}\|_{L^2}^{\delta_4}\right)\left(\lambda_p^{s(2-\delta_4)}\| B_{p}\|_{L^2}^{2-\delta_4}\right)  \\
&\cdot \sum_{q\leq p+2}c_q\lambda_q^{s}\|B_q\|_{L^2} \lambda_{p-q}^{1+r_1+r_2-2s-\frac\alpha2\delta_4}\\
\leq&\ \frac{1}{64}\mu\sum_{p\geq -1}\lambda_p^{2s+\alpha}\| B_{p}\|_{L^2}^{2}\\
&+\frac{C}{\mu}\sum_{p\geq -1}\lambda_p^{2s}\| B_{p}\|_{L^2}^{2}\left(  \sum_{q\leq p+2}c_q\lambda_q^{s}\|B_q\|_{L^2}\lambda_{p-q}^{1+r_1+r_2-2s-\frac\alpha2\delta_4}\right)^{\frac2{2-\delta_4}}\\
\leq&\ \frac{1}{64}\mu\sum_{p\geq -1}\lambda_p^{2s+\alpha}\| B_{p}\|_{L^2}^{2}+\frac{C}{\mu}\left(\sum_{p\geq -1}\lambda_p^{2s}\| B_{p}\|_{L^2}^{2}\right)^{1+\frac1{2-\delta_4}}.
\end{split}
\end{equation}
Summarizing the analysis above, we can choose $\delta_1=\delta_2=\delta_3=\delta_4=\delta\in [\frac{5-2s}{\alpha},2)$ and close enough to 2, and conclude 
\begin{equation}\label{est-I}
|I|
\leq\frac{1}{4}\mu\sum_{q\geq -1}\lambda_q^{2s+\alpha}\| B_{q}\|_{L^2}^{2}+\frac{C}{\mu}\left(\sum_{q\geq -1}\lambda_q^{2s}\| B_{q}\|_{L^2}^{2}\right)^{1+\frac1{2-\delta}}.
\end{equation}

 \subsection{Estimate of $K$}

Applying Bony's paraproduct we decompose $K$ into
\begin{equation}\notag
\begin{split}
K=&\sum_{q\geq -1}\sum_{|p-q|\leq 2}\lambda_q^{2s}\int_{\mathbb R}\Delta_q(\Lambda B_{p}B_{x,\leq p-2})B_q dx\\
&+\sum_{q\geq -1}\sum_{|p-q|\leq 2}\lambda_q^{2s}\int_{\mathbb R}\Delta_q(\Lambda B_{\leq p-2}B_{x,p})B_q dx\\
&+\sum_{q\geq -1}\sum_{p\geq q-2}\lambda_q^{2s}\int_{\mathbb R}\Delta_q(\Lambda \widetilde B_{p}B_{x,p})B_q dx\\
=&: K_{1}+K_{2}+K_{3}.
\end{split}
\end{equation}
It follows from H\"older's inequality and Bernstein's inequality that
\begin{equation}\notag
\begin{split}
|K_{1}|\lesssim & \sum_{q\geq -1}\lambda_q^{2s}\|B_q\|_{L^2} \sum_{|p-q|\leq 2}\|\Lambda B_{p}\|_{L^2} \|B_{x, \leq p-2}\|_{L^\infty} \\
\lesssim & \sum_{q\geq -1}\lambda_q^{2s+1}\|B_q\|_{L^2}^2 \sum_{p\leq q+2}\lambda_p^{\frac32}\|B_{p}\|_{L^2} \\
\end{split}
\end{equation}
which has a similar form as for $I_{11}$ in \eqref{est-I11} and hence has a similar estimate as in \eqref{est-I11-b},
\begin{equation}\notag
|K_{1}|
\leq \frac{1}{64}\mu\sum_{q\geq -1}\lambda_q^{2s+\alpha}\|B_q\|_{L^2}^2+\frac{C}{\mu}\left(\sum_{q\geq -1}\lambda_q^{2s}\|B_q\|_{L^2}^2\right)^{1+\frac1{2-\delta_1}}.
\end{equation}
The term $K_{2}$ can be estimated analogously as for $K_{1}$. 

To estimate $K_{3}$, we need to further decompose it applying a commutator as
\begin{equation}\notag
\begin{split}
K_{3}=& \sum_{q\geq -1}\sum_{p\geq q-2}\lambda_q^{2s}\int_{\mathbb R}\Lambda \widetilde B_{p}\Delta_qB_{x,p}B_q dx\\
&+\sum_{q\geq -1}\sum_{p\geq q-2}\lambda_q^{2s}\int_{\mathbb R}[\Delta_q, \Lambda \widetilde B_{p}]B_{x,p}B_q dx\\
=&: K_{31}+K_{32}.
\end{split}
\end{equation}
%where the commutator is defined as
%\[[\Delta_q, \Lambda \widetilde B_{p}]B_{x,p}=\Delta_q(\Lambda \widetilde B_{p}B_{x,p})-\Lambda \widetilde B_{p}\Delta_qB_{x,p}.\]
Observing that $\Delta_qB_{x,p}=0$ for $|p-q|\geq 3$, the summation $\sum_{p\geq q-2}$ has only a few non-zero terms. Applying H\"older's inequality and Bernstein's inequality again we have
\begin{equation}\notag
\begin{split}
|K_{31}|\lesssim &\sum_{q\geq -1}\sum_{p\sim q}\lambda_q^{2s}\|\Lambda \widetilde B_p\|_{L^2}\|\Delta_qB_{x,p}\|_{L^2}\|B_q\|_{L^\infty}\\
\lesssim &\sum_{q\geq -1}\lambda_q^{2s+\frac52}\|B_q\|_{L^2}^3\\
\leq &\ \frac{1}{64}\mu\sum_{q\geq -1}\lambda_q^{2s+\alpha}\|B_q\|_{L^2}^2+\frac{C}{\mu}\left( \sum_{q\geq -1}\lambda_q^{2s}\|B_q\|_{L^2}^2\right)^{\frac{3-\delta_3}{2-\delta_3}}\\
\end{split}
\end{equation}
which is obtained similarly as for $I_{13}$ as in \eqref{est-I13-b}.

To estimate $K_{32}$, we apply H\"older's inequality and Lemma \ref{le-1} with $r_0=0$, $n=1$, $\frac12\leq r_1<\frac32$, $r_2<\frac12$ and $r_1+r_2>0$
\begin{equation}\notag
\begin{split}
|K_{32}|\lesssim & \sum_{q\geq -1}\lambda_q^{2s}\|B_q\|_{L^2}\sum_{p\geq q-2}\|[\Delta_q, \Lambda \widetilde B_{p}]B_{x,p}\|_{L^2}\\
\lesssim & \sum_{q\geq -1}c_q\lambda_q^{2s-(r_1+r_2-\frac12)}\|B_q\|_{L^2}\sum_{p\geq q-2}\|\Lambda \widetilde B_{p}\|_{\dot H^{r_1}}\|B_{x,p}\|_{\dot H^{r_2}}\\
\lesssim & \sum_{p\geq -1}\|\Lambda \widetilde B_{p}\|_{\dot H^{r_1}}\|B_{x,p}\|_{\dot H^{r_2}}\sum_{q\leq p+2}c_q\lambda_q^{2s-(r_1+r_2-\frac12)}\|B_q\|_{L^2}\\
\lesssim & \sum_{p\geq -1}\lambda_p^{2+r_1+r_2}\| B_{p}\|_{L^2}^2\sum_{q\leq p+2}c_q\lambda_q^{2s-(r_1+r_2-\frac12)}\|B_q\|_{L^2}.
\end{split}
\end{equation}
We rearrange the terms on the right hand side for some $\delta_5\in(0,2)$  
\begin{equation}\notag
\begin{split}
|K_{32}|
\lesssim & \sum_{p\geq -1}\lambda_p^{2+r_1+r_2}\| B_{p}\|_{L^2}^2\sum_{q\leq p+2}c_q\lambda_q^{2s-(r_1+r_2-\frac12)}\|B_q\|_{L^2}\\
\lesssim & \sum_{p\geq -1}\left(\lambda_p^{(s+\frac\alpha2)\delta_5}\| B_{p}\|_{L^2}^{\delta_5}\right)\left(\lambda_p^{s(2-\delta_5)}\| B_{p}\|_{L^2}^{2-\delta_5}\right)  \sum_{q\leq p+2}c_q\lambda_q^{s}\|B_q\|_{L^2}\\
&\cdot \lambda_q^{\frac52-s-\frac\alpha2\delta_5}\lambda_{p-q}^{2+r_1+r_2-2s-\frac\alpha2\delta_5}.
\end{split}
\end{equation}
To optimize the estimate, we choose $r_1=\frac12$ and $r_2=-\frac12+\varepsilon$ for a small constant $\varepsilon>0$ such that
\begin{equation}\label{para-5}
\begin{split}
\frac52-s-\frac\alpha2\delta_5&\leq 0,\\
2+r_1+r_2-2s-\frac\alpha2\delta_5&<0.
\end{split}
\end{equation}
Note for $s>\frac52-\alpha$ and $\alpha\in(0,\frac52]$, we can choose $\delta_5\in[\frac{5-2s}{\alpha}, 2)$ and $\delta_5$ close enough to 2 to have the conditions in \eqref{para-5} fulfilled. 
Therefore we proceed the estimate by applying Young's inequality and Jensen's inequality 
\begin{equation}\notag
\begin{split}
|K_{32}|
\lesssim & \sum_{p\geq -1}\left(\lambda_p^{(s+\frac\alpha2)\delta_5}\| B_{p}\|_{L^2}^{\delta_3}\right)\left(\lambda_p^{s(2-\delta_5)}\| B_{p}\|_{L^2}^{2-\delta_5}\right)  \sum_{q\leq p+2}c_q\lambda_q^{s}\|B_q\|_{L^2}\\
&\cdot \lambda_q^{\frac52-s-\frac\alpha2\delta_5}\lambda_{p-q}^{-\frac\alpha2\delta_5}\\
\lesssim & \sum_{p\geq -1}\left(\lambda_p^{(s+\frac\alpha2)\delta_5}\| B_{p}\|_{L^2}^{\delta_5}\right)\left(\lambda_p^{s(2-\delta_5)}\| B_{p}\|_{L^2}^{2-\delta_5}\right)  \sum_{q\leq p+2}c_q\lambda_q^{s}\|B_q\|_{L^2}\lambda_{p-q}^{-\frac\alpha2\delta_5}\\
\leq&\ \frac{1}{64}\mu\sum_{p\geq -1}\lambda_p^{2s+\alpha}\| B_{p}\|_{L^2}^{2}+\frac{C}{\mu}\sum_{p\geq -1}\lambda_p^{2s}\| B_{p}\|_{L^2}^{2}\left(  \sum_{q\leq p+2}c_q\lambda_q^{s}\|B_q\|_{L^2}\lambda_{p-q}^{-\frac\alpha2\delta_5}\right)^{\frac2{2-\delta_5}}\\
\leq&\ \frac{1}{64}\mu\sum_{p\geq -1}\lambda_p^{2s+\alpha}\| B_{p}\|_{L^2}^{2}+\frac{C}{\mu}\left(\sum_{p\geq -1}\lambda_p^{2s}\| B_{p}\|_{L^2}^{2}\right)^{1+\frac1{2-\delta_5}}.
\end{split}
\end{equation}

Combining the estimates above and choosing $\delta_1=\delta_3=\delta_5=\delta\in[\frac{5-2s}{\alpha}, 2)$ which is close enough to 2, we obtain
\begin{equation}\label{est-K}
|K|\leq \frac18\mu \sum_{q\geq -1}\lambda_q^{2s+\alpha}\| B_{q}\|_{L^2}^{2}+\frac{C}{\mu}\left(\sum_{q\geq -1}\lambda_q^{2s}\| B_{q}\|_{L^2}^{2}\right)^{1+\frac1{2-\delta}}.
\end{equation}

In the end, in view of \eqref{est-energy}, \eqref{est-I} and \eqref{est-K}, we claim
\begin{equation}\label{est-energy2}
\frac{d}{dt}\sum_{q\geq -1}\lambda_q^{2s}\|B_q\|_{L^2}^2+\mu\sum_{q\geq -1}\lambda_q^{2s+\alpha} \|B_q\|_{L^2}^2
\leq  \frac{C}{\mu}\left(\sum_{q\geq -1}\lambda_q^{2s}\| B_{q}\|_{L^2}^{2}\right)^{1+\frac1{2-\delta}}.
\end{equation}
It then follows from \eqref{est-energy2} that there exists a time $T>0$ depending on the $\dot H^s$ norm of the initial profile $B(x,0)$, $\delta$ and $\mu$ such that
\begin{equation}\notag
\sum_{q\geq -1}\lambda_q^{2s}\|B_q(t)\|_{L^2}^2\leq \frac{\sum_{q\geq -1}\lambda_q^{2s}\|B_q(0)\|_{L^2}^2}{\left(1-C(2-\delta)^{-1}\mu^{-1}t (\sum_{q\geq -1}\lambda_q^{2s}\|B_q(0)\|_{L^2}^2)^{\frac1{2-\delta}} \right)^{2-\delta}}
\end{equation}
for $0\leq t<T$, and 
\[\int_0^t\sum_{q\geq -1}\lambda_q^{2s+\alpha} \|B_q(\tau)\|_{L^2}^2\, d\tau\leq C_0\sum_{q\geq -1}\lambda_q^{2s}\|B_q(0)\|_{L^2}^2 \]
with the constant $C_0$ depending on $\|B(0)\|_{\dot H^s}$, $\delta$ and $\mu$. Hence we obtain the a priori estimate
\begin{equation}\label{priori}
B\in C\left([0,T);   \dot H^{s}(\mathbb R)\right)\cap L^2\left([0,T);  \dot H^{s+\frac{\alpha}2}(\mathbb R)\right), \ s>\frac52-\alpha.
\end{equation}

\medskip

\subsection{A priori estimate in $L^2$} 
\label{sec-l2}

%Question: {\color{red} We are on torus, do we need this part below to show estimate in $L^2$? On torus, does estimate in $H^s$ with $s>0$ implies estimate in $L^2$?}
We now apply the higher order a priori estimate \eqref{priori} established in previous section to obtain the a priori estimate in $L^2$.
Taking $s=0$ in (\ref{est-energy}) gives
\begin{equation}\label{est-l2}
\begin{split}
&\frac12\frac{d}{dt}\sum_{q\geq -1}\|B_q\|_{L^2}^2+\mu\sum_{q\geq -1} \|\Lambda^{\frac\alpha2}B_q\|_{L^2}^2\\
=&-\sum_{q\geq -1}\int_{\mathbb R}(B\Lambda B)_qB_{x,q} dx-2\sum_{q\geq -1}\int_{\mathbb R}(\Lambda BB_x)_qB_q dx\\
=&: -\bar I-2 \bar K
\end{split}
\end{equation}
with Bony's paraproduct of $\bar I$ and $\bar K$ given by
\begin{equation}\notag
\begin{split}
\bar I=&\sum_{q\geq -1}\sum_{|p-q|\leq 2}\int_{\mathbb R}\Delta_q(B_{\leq p-2}\Lambda B_{p})B_{x,q} dx\\
&+\sum_{q\geq -1}\sum_{|p-q|\leq 2}\int_{\mathbb R}\Delta_q(B_{p}\Lambda B_{\leq p-2})B_{x,q} dx\\
&+\sum_{q\geq -1}\sum_{p\geq q-2}\int_{\mathbb R}\Delta_q(B_{p}\Lambda \widetilde B_{p})B_{x,q} dx\\
=&:\bar I_{1}+\bar I_{2}+\bar I_{3},
\end{split}
\end{equation}
and
\begin{equation}\notag
\begin{split}
\bar K=&\sum_{q\geq -1}\sum_{|p-q|\leq 2}\int_{\mathbb R}\Delta_q(\Lambda B_{p}B_{x,\leq p-2})B_q dx\\
&+\sum_{q\geq -1}\sum_{|p-q|\leq 2}\int_{\mathbb R}\Delta_q(\Lambda B_{\leq p-2}B_{x,p})B_q dx\\
&+\sum_{q\geq -1}\sum_{p\geq q-2}\int_{\mathbb R}\Delta_q(\Lambda \widetilde B_{p}B_{x,p})B_q dx\\
=&: \bar K_{1}+\bar K_{2}+\bar K_{3}.
\end{split}
\end{equation}

Applying H\"older's inequality and Bernstein's inequality yields
\begin{equation}\label{bar-I1}
\begin{split}
|\bar I_1|\leq& \sum_{q\geq -1}\sum_{|p-q|\leq 2}\|\Lambda B_p\|_{L^2}\|B_{x,q}\|_{L^2}\|B_{\leq p-2}\|_{L^\infty}\\
\lesssim& \sum_{q\geq -1}\lambda_q^2\|B_{q}\|_{L^2}^2\sum_{p\leq q}\lambda_p^{\frac12}\|B_{p}\|_{L^2}.
\end{split}
\end{equation}
We then apply Young's inequality and Jensen's inequality for $s\geq \frac52-\frac\alpha2$ to infer
\begin{equation}\notag
\begin{split}
|\bar I_1|
\lesssim& \sum_{q\geq -1}\left(\lambda_q^{\frac\alpha2}\|B_{q}\|_{L^2}\right)\left(\lambda_q^s\|B_{q}\|_{L^2}\right) \sum_{p\leq q}\|B_{p}\|_{L^2}\lambda_{p-q}^{\frac12}\lambda_q^{\frac52-\frac\alpha2-s}\\
\leq&\ \frac{1}{64}\sum_{q\geq -1}\lambda_q^{\alpha}\|B_{q}\|_{L^2}^2+\frac{C}\mu \sum_{q\geq -1}\lambda_q^{2s}\|B_{q}\|_{L^2}^2\left( \sum_{p\leq q}\|B_{p}\|_{L^2}\lambda_{p-q}^{\frac12}\right)^2\\
\leq&\ \frac{1}{64}\sum_{q\geq -1}\lambda_q^{\alpha}\|B_{q}\|_{L^2}^2+\frac{C}\mu \sum_{q\geq -1}\lambda_q^{2s}\|B_{q}\|_{L^2}^2 \sum_{p\leq q}\|B_{p}\|_{L^2}^2\\
\leq&\ \frac{1}{64}\sum_{q\geq -1}\lambda_q^{\alpha}\|B_{q}\|_{L^2}^2+\frac{C}\mu  \sum_{p\geq -1}\|B_{p}\|_{L^2}^2\sum_{q\geq p}\lambda_q^{2s}\|B_{q}\|_{L^2}^2\\
\leq&\ \frac{1}{64}\sum_{q\geq -1}\lambda_q^{\alpha}\|B_{q}\|_{L^2}^2+\frac{C}\mu \|B\|_{\dot H^s}^2 \sum_{p\geq -1}\|B_{p}\|_{L^2}^2.
\end{split}
\end{equation}
Similarly we estimate $\bar I_2$ for $s\geq \frac52-\frac\alpha2$
\begin{equation}\notag
\begin{split}
|\bar I_2|\leq& \sum_{q\geq -1}\sum_{|p-q|\leq 2}\| B_p\|_{L^2}\|B_{x,q}\|_{L^2}\|\Lambda B_{\leq p-2}\|_{L^\infty}\\
\lesssim& \sum_{q\geq -1}\lambda_q\|B_{q}\|_{L^2}^2\sum_{p\leq q}\lambda_p^{\frac32}\|B_{p}\|_{L^2}\\
\lesssim& \sum_{q\geq -1}\left(\lambda_q^{\frac\alpha2}\|B_{q}\|_{L^2}\right)\left(\lambda_q^s\|B_{q}\|_{L^2}\right) \sum_{p\leq q}\|B_{p}\|_{L^2}\lambda_{p-q}^{\frac32}\lambda_q^{\frac52-\frac\alpha2-s}\\
%\leq&\ \frac{1}{64}\sum_{q\geq -1}\lambda_q^{\alpha}\|B_{q}\|_{L^2}^2+\frac{C}\mu \sum_{q\geq -1}\lambda_q^{2s}\|B_{q}\|_{L^2}^2\left( \sum_{p\leq q}\|B_{p}\|_{L^2}\lambda_{p-q}^{\frac12}\right)^2\\
%\leq&\ \frac{1}{64}\sum_{q\geq -1}\lambda_q^{\alpha}\|B_{q}\|_{L^2}^2+\frac{C}\mu \sum_{q\geq -1}\lambda_q^{2s}\|B_{q}\|_{L^2}^2 \sum_{p\leq q}\|B_{p}\|_{L^2}^2\\
%\leq&\ \frac{1}{64}\sum_{q\geq -1}\lambda_q^{\alpha}\|B_{q}\|_{L^2}^2+\frac{C}\mu  \sum_{p\geq -1}\|B_{p}\|_{L^2}^2\sum_{q\geq p}\lambda_q^{2s}\|B_{q}\|_{L^2}^2\\
\leq&\ \frac{1}{64}\sum_{q\geq -1}\lambda_q^{\alpha}\|B_{q}\|_{L^2}^2+\frac{C}\mu \|B\|_{\dot H^s}^2 \sum_{p\geq -1}\|B_{p}\|_{L^2}^2.
\end{split}
\end{equation}
We observe $\bar I_3$, $\bar K_1$ and $\bar K_2$ share similar estimate of $\bar I_2$. Indeed, we have
\begin{equation}\notag
\begin{split}
|\bar I_3|\leq& \sum_{q\geq -1}\sum_{p\geq q-2}\|B_{p}\|_{L^2}\|\Lambda \widetilde B_p\|_{L^2}\|B_{x,q}\|_{L^\infty}\\
\lesssim& \sum_{q\geq -1}\lambda_q^{\frac32}\|B_{q}\|_{L^2}\sum_{p\geq q-2}\lambda_p\|B_{p}\|_{L^2}^2\\
\lesssim& \sum_{p\geq -1}\lambda_p\|B_{p}\|_{L^2}^2\sum_{q\leq p+2}\lambda_q^{\frac32}\|B_{q}\|_{L^2}\\
\leq&\ \frac{1}{64}\sum_{q\geq -1}\lambda_q^{\alpha}\|B_{q}\|_{L^2}^2+\frac{C}\mu \|B\|_{\dot H^s}^2 \sum_{p\geq -1}\|B_{p}\|_{L^2}^2,
\end{split}
\end{equation}
and
\begin{equation}\notag
\begin{split}
|\bar K_1|+|\bar K_2|\leq& \sum_{q\geq -1}\sum_{|p-q|\leq 2}\|\Lambda B_p\|_{L^2}\|B_{q}\|_{L^2}\|B_{x,\leq p-2}\|_{L^\infty}\\
&+\sum_{q\geq -1}\sum_{|p-q|\leq 2}\| B_q\|_{L^2}\|B_{x,p}\|_{L^2}\|\Lambda B_{\leq p-2}\|_{L^\infty}\\
\lesssim& \sum_{q\geq -1}\lambda_q\|B_{q}\|_{L^2}^2\sum_{p\leq q}\lambda_p^{\frac32}\|B_{p}\|_{L^2}\\
\leq&\ \frac{1}{64}\sum_{q\geq -1}\lambda_q^{\alpha}\|B_{q}\|_{L^2}^2+\frac{C}\mu \|B\|_{\dot H^s}^2 \sum_{p\geq -1}\|B_{p}\|_{L^2}^2.
\end{split}
\end{equation}
%\begin{equation}\notag
%\begin{split}
%|\bar K_2|\leq& \sum_{q\geq -1}\sum_{|p-q|\leq 2}\| B_q\|_{L^2}\|B_{x,p}\|_{L^2}\|\Lambda B_{\leq p-2}\|_{L^\infty}\\
%\lesssim& \sum_{q\geq -1}\lambda_q\|B_{q}\|_{L^2}^2\sum_{p\leq q}\lambda_p^{\frac32}\|B_{p}\|_{L^2}\\
%\leq&\ \frac{1}{64}\sum_{q\geq -1}\lambda_q^{\alpha}\|B_{q}\|_{L^2}^2+\frac{C}\mu \|B\|_{\dot H^s}^2 \sum_{p\geq -1}\|B_{p}\|_{L^2}^2\\
%\end{split}
%\end{equation}
Finally we note $\bar K_3$ can be estimated analogously as for $\bar I_1$, in view of \eqref{bar-I1}
\begin{equation}\notag
\begin{split}
|\bar K_3|\leq& \sum_{q\geq -1}\sum_{p\geq q-2}\|B_{x,p}\|_{L^2}\|\Lambda \widetilde B_p\|_{L^2}\|B_{q}\|_{L^\infty}\\
\lesssim& \sum_{q\geq -1}\lambda_q^{\frac12}\|B_{q}\|_{L^2}\sum_{p\geq q-2}\lambda_p^2\|B_{p}\|_{L^2}^2\\
\lesssim& \sum_{p\geq -1}\lambda_p^2\|B_{p}\|_{L^2}^2\sum_{q\leq p+2}\lambda_q^{\frac12}\|B_{q}\|_{L^2}\\
\leq&\ \frac{1}{64}\sum_{q\geq -1}\lambda_q^{\alpha}\|B_{q}\|_{L^2}^2+\frac{C}\mu \|B\|_{\dot H^s}^2 \sum_{p\geq -1}\|B_{p}\|_{L^2}^2.
\end{split}
\end{equation}
Combining the estimates of $\bar I$ and $\bar K$ together with \eqref{est-l2} gives
\begin{equation}\notag%\label{est-l2-final}
\frac{d}{dt}\sum_{q\geq -1}\|B_q\|_{L^2}^2+\mu\sum_{q\geq -1}\lambda_q^{\alpha} \|B_q\|_{L^2}^2
\leq \frac{C}\mu \|B\|_{\dot H^s}^2 \sum_{p\geq -1}\|B_{p}\|_{L^2}^2.
\end{equation}
It then follows from \eqref{priori} for $s>\frac52-\frac\alpha2$ that 
\begin{equation}\label{priori2}
B\in C\left([0,T_1);  L^{2}(\mathbb R)\right)\cap L^2\left([0,T_1);  \dot H^{\frac{\alpha}2}(\mathbb R)\right)
\end{equation}
for some $0<T_1<T$.

%The estimates \eqref{priori} and \eqref{priori2} together justify \eqref{priori0}.

\medskip

\subsection{Proof of Theorem \ref{thm-local}}
\label{sec-proof1}
To show the existence of the solution stated in Theorem \ref{thm-local}, we follow the standard procedure as in \cite{Dai-1d-emhd}. Namely
we solve the approximating system iteratively
\begin{equation}\label{sys-app}
\begin{split}
B^{k}_t+B^{k-1} J_x^{k}-J^{k-1}B_x^{k}+\mu \Lambda^\alpha B^{k}=&\ 0,\\
B_x^k=&\ \mathcal H J^k,\\
B^{k}(x,0)=&\ B_0(x)
\end{split}
\end{equation}
for $k\geq 0$. By convention, take $B^{-1}=J^{-1}\equiv 0$. The sequence of solutions $\{B^k\}$ can be shown to satisfy the a priori estimates \eqref{priori} and \eqref{priori2}. 
Consequently, one can show that a subsequence of $\{B^k\}$ converges to a limit function $B$ in 
$L^\infty([0, T); H^{s})\cap L^2([0,T); H^{s+\frac{\alpha}{2}})$ for $s>\frac52-\alpha$
and $B$ is a classical solution of (\ref{emhd-1d}) on $[0,T)$. The continuity in time and smoothing estimates \eqref{est-b1}-\eqref{est-b2} can be obtained analogously as in \cite{Dai-1d-emhd}.

\cbdu

%\bigskip

%\section{Analytic solutions for small $\alpha$}

%\bigskip

%\section{Local smooth solutions}
%Consider the ``viscosity" solutions first; use weak compactness to take the limit to obtain solutions for the ``inviscid" case. Refer to Tao's notes.

\bigskip

\section{Blowup in the resistive case without stretching effect}

This section provides a proof for Theorem \ref{thm-blow}.

%{\textbf{Proof of Theorem \ref{thm-blow}:}}
%Without loss of generality, we take the resistivity parameter $\mu=1$. 
Recall the equation \eqref{emhd-transport} with $\mu=1$ and $\alpha=1$
%\[\partial_tB+JB_x+\Lambda B=0.\]
%Since $J=-\mathcal H\partial_xB=-\Lambda B$, we have
\begin{equation}\label{B-2}
\partial_tB-\Lambda BB_x+\Lambda B=0.
\end{equation}
Denote $\bar B(x,t)=B(x,t)-x$ which satisfies the transport equation
\begin{equation}\label{B-bar}
\partial_t\bar B-\Lambda B\bar B_x=0.
\end{equation}
Let $X(x_0, t)$ be the trajectory defined by
\begin{equation}\label{traj1}
\begin{split}
\frac{d}{dt} X(x_0, t)%=J(X(x_0, t), t)
&=-\Lambda B(X(x_0, t), t),\\
X(x_0, 0)&=x_0.
\end{split}
\end{equation}
Applying \eqref{B-bar} and \eqref{traj1} gives
\begin{equation}\notag
\begin{split}
&\frac{d}{dt}\bar B_x(X(x_0,t),t)\\
=&\ \partial_t \bar B_x(X(x_0,t),t)+\frac{d}{dt} X(x_0, t)\bar B_{xx}(X(x_0,t),t)\\
=&\ \Lambda B_x(X(x_0,t),t)\bar B_x(X(x_0,t),t)+\Lambda B(X(x_0,t),t)\bar B_{xx}(X(x_0,t),t)\\
& -\Lambda B(X(x_0,t),t)\bar B_{xx}(X(x_0,t),t),
\end{split}
\end{equation}
which is
\begin{equation}\label{B-bar-x}
\frac{d}{dt}\bar B_x(X(x_0,t),t)=\Lambda B_x(X(x_0,t),t)\bar B_x(X(x_0,t),t).
\end{equation}
Applying \eqref{B-2} and \eqref{traj1} we can compute
\begin{equation}\notag
\begin{split}
&\frac{d}{dt}\mathcal H B_{xx}(X(x_0,t),t)\\
=&\ \partial_t \mathcal H B_{xx}(X(x_0,t),t)+\frac{d}{dt} X(x_0, t)\mathcal H B_{xxx}(X(x_0,t),t)\\
=&\ \mathcal H\partial_{xx}\left(\Lambda B(X(x_0,t),t)B_x(X(x_0,t),t)-\Lambda B(X(x_0,t),t) \right)\\
&-\Lambda B(X(x_0,t),t)\mathcal H B_{xxx}(X(x_0,t),t)\\
=&\ \partial_{xx}\mathcal H\left(\mathcal H B_x(X(x_0,t),t)B_x(X(x_0,t),t)) \right)\\
&+ B_{xxx}(X(x_0,t),t) -\Lambda B(X(x_0,t),t)\mathcal H B_{xxx}(X(x_0,t),t)
\end{split}
\end{equation}
where we used $\Lambda=\mathcal H\partial_x$. 
In view of the property of Hilbert transform
\begin{equation}\label{fhf}
\mathcal H(f\mathcal H f)=\frac12\left[(\mathcal H f)^2-f^2 \right],
\end{equation}
we continue to infer
\begin{equation}\notag
\begin{split}
&\frac{d}{dt}\mathcal H B_{xx}(X(x_0,t),t)\\
=&\ \partial_{xx}\left(\frac12(\mathcal H B_x(X(x_0,t),t))^2-\frac12(B_x(X(x_0,t),t))^2 \right)\\
&+ B_{xxx}(X(x_0,t),t) -\Lambda B(X(x_0,t),t)\mathcal H B_{xxx}(X(x_0,t),t)\\
=&\ \left(\mathcal HB_{xx}(X(x_0,t),t)\right)^2+\mathcal HB_x(X(x_0,t),t)\mathcal HB_{xxx}(X(x_0,t),t)\\
&-\left(B_{xx}(X(x_0,t),t)\right)^2-B_x(X(x_0,t),t)B_{xxx}(X(x_0,t),t)\\
&+ B_{xxx}(X(x_0,t),t) -\Lambda B(X(x_0,t),t)\mathcal H B_{xxx}(X(x_0,t),t).
\end{split}
\end{equation}

Hence we have
\begin{equation}\label{B-xx}
\begin{split}
&\frac{d}{dt}\mathcal H B_{xx}(X(x_0,t),t)\\
=&\ \left(\mathcal HB_{xx}(X(x_0,t),t)\right)^2-\left(B_{xx}(X(x_0,t),t)\right)^2\\
&-\left(B_x(X(x_0,t),t)-1\right)B_{xxx}(X(x_0,t),t).
\end{split}
\end{equation}
We observe that, if the second and third terms on the right hand side of \eqref{B-xx} vanish along the trajectory $X(x_0,t)$, the equation is in Riccati type and the solution of which would blow up at a finite time. Therefore, it motivates us to choose the initial profile $B_0(x)$ such that $\partial_x B_0(x_0)=1$, $\partial_{xx} B_0(x_0)=0$ and $\mathcal H\partial_{xx}B_0(x_0)>0$ for some $x_0$. Since $\partial_x\bar B(x_0)=0$, equation \eqref{B-bar-x} implies 
\[\partial_x\bar B(X(x_0,t),t)=0 \ \ \mbox{for all} \ \ t\geq 0 \]
and hence 
\[\partial_x B(X(x_0,t),t)=1 \ \ \mbox{for all} \ \ t\geq 0.\]
Thus it follows from \eqref{B-xx} that
\begin{equation}\notag
\frac{d}{dt}\mathcal H B_{xx}(X(x_0,t),t)
=\left(\mathcal HB_{xx}(X(x_0,t),t)\right)^2-\left(B_{xx}(X(x_0,t),t)\right)^2.
\end{equation}
Taking the Hilbert transform on the equation above and applying \eqref{fhf} again, we obtain
\begin{equation}\label{B-xx2}
\frac{d}{dt} B_{xx}(X(x_0,t),t)
=2B_{xx}(X(x_0,t),t)\mathcal HB_{xx}(X(x_0,t),t).
\end{equation}
Since $\partial_{xx} B_0(x_0)=0$, it follows from equation \eqref{B-xx2} that
\[\partial_{xx} B(X(x_0,t),t)=0 \ \ \mbox{for all} \ \ t\geq 0.\]
Therefore equation \eqref{B-xx} is essentially a Riccati type of equation
\begin{equation}\notag
\frac{d}{dt}\mathcal H B_{xx}(X(x_0,t),t)
= \left(\mathcal HB_{xx}(X(x_0,t),t)\right)^2
\end{equation}
whose solution is given by
\[\mathcal H B_{xx}(X(x_0,t),t)=\frac{\mathcal H \partial_{xx}B_0(x_0)}{1-t \mathcal H \partial_{xx}B_0(x_0)}.\]
Thus $\mathcal H B_{xx}$ blows up along the trajectory $X(x_0,t)$ at the time $T=\frac1{\mathcal H \partial_{xx}B_0(x_0)}>0$. It completes the proof.

\cbdu

\bigskip

%\section{Pure transport case}
%\label{sec-transport}

%\section*{Acknowledgement}
%The author is indebted to Diego C\'ordoba and Hongjie Dong for their valuable suggestions which have improved significantly the early version of the article. 

%The author would like to express her sincere gratitude towards the anonymous referee who has read the first manuscript carefully and provided constructive suggestions. 

%\Endrefs
\end{document}